\documentclass[a4paper,fleqn]{cas-dc}

\usepackage[numbers]{natbib}

\def\tsc#1{\csdef{#1}{\textsc{\lowercase{#1}}\xspace}}
\tsc{WGM}
\tsc{QE}
\tsc{EP}
\tsc{PMS}
\tsc{BEC}
\tsc{DE}
\newcommand{\vect}[1]{\boldsymbol{\mathit{#1}}}
\newcommand{\tenf}[1]{\boldsymbol{\mathbb{#1}}}%

\begin{document}
\let\WriteBookmarks\relax
\def\floatpagepagefraction{1}
\def\textpagefraction{.001}
\shorttitle{Non-Hookean elasticity with arbitrary Poisson's ratios}
\shortauthors{M. Itskov}

\title [mode = title]{Non-Hookean elasticity with arbitrary Poisson's ratios}

\author[1]{Mikhail Itskov}[type=editor,
                        auid=000,bioid=1,
                        orcid=0000-0002-2146-0589]
\cormark[1]

\ead{itskov@km.rwth-aachen.de}
\ead[url]{www.km.rwth-aachen.de}

\affiliation[1]{organization={Department of Continuum Mechanics, RWTH Aachen University},
                addressline={Eilfschornsteinstr. 18}, 
                city={Aachen},
                postcode={52062}, 
                country={Germany}}

\cortext[cor1]{Corresponding author}

\begin{abstract}
In a previous paper \cite{Itskov-MoSM} we presented a hyperelastic isotropic material model whose stress-strain response is nonlinear even at infinitesimal deformations and cannot thus be linearized. As a result values of Poisson's ratio greater than one half were obtained. In this contribution, we further propose an isotropic strain energy function which is always positive-definite and depending on material constants delivers arbitrary values of Poisson's ratio (except of $-1$) in agreement with the laws of thermodynamics. The model response appears stable and plausible in various deformation states.
\end{abstract}

\begin{highlights}
\item Non-linearizable hyperelastic isotropic model based on a positive-definite strain energy function is proposed
\item Non-Hookean elasticity theory without superposition principle is presented
\item Poisson's ratios beyond the classical range $[-1,0.5]$ for isotropic materials are predicted in agreement with the thermodynamics laws 
\end{highlights}


\maketitle

\section{Introduction}

{Hookean elasticity implies a linear relation between stresses and strains, which is reasonable for most of solids at infinitesimal deformations. In turn, linearity justifies the superposition principle leading to the classical generalized Hooke law (Lam\'{e} formulas). One of its well-known thermodynamical consequences is that  for isotropic materials Poisson's 
	ratio $\nu$ lies within the interval  $\left[-1,0.5\right] $ (see, e.g. \cite{True}, Sect. 51).  For example, values of $\nu$ smaller than $-1$ would lead to the negative shear modulus and infinite release of energy under simple shear (see e.g. \cite{Landau,Love,Mott}). 
	
	In this contribution, we propose a hyperelastic isotropic material model whose stress-strain response is nonlinear even at infinitesimal deformations and cannot thus be linearized.
	As a result, the superposition principle does not apply and the generalized Hooke law is not valid even at small strains. In this case, one can speak about non-Hookean elasticity both at small and finite strains.
	A further unusual feature of this material model is that Poisson's ratio can be greater than $0.5$ or smaller than $-1$
	in full agreement with the laws of thermodynamics. Thus, depending
	on material constants arbitrary values of Poisson's ratio except of $-1$ are predicted by the model. We also discuss the issue with the value of $-1$ and demonstrate that it is only realizable in an anisotropic material.
	
	In a series of papers (see, e.g. \cite{Raja-impl}) Rajagopal proposed implicit constitutive relations of the form $f\left(\tenf{\sigma},\tenf{\epsilon}\right)=0$, where $\tenf{\sigma}$ and $\tenf{\epsilon}$ denote the Cauchy stress and strain tensor, respectively. Recently, this concept was mathematically justified, see \cite{Raja-recent}. In the references of the latter paper one can also find many examples of isotropic materials whose nonlinear response at small strains was confirmed in experiments.
	
	Within this concept Rajagopal considered in \cite{Raja} 
	a strain response as a quadratic function of stress. Although this non-linearity persists at small strains it would disappear at small stresses. Thus, under linearization both with respect to strains and stresses his formulation will still imply the superposition principle and reduce to generalized Hooke's law.
	
	The paper is organized as follows. In Sect. \ref{Model} we present the general concept and formulate a positive-definite strain energy function. This strain energy function leads to a constitutive relation similar
	to the force-displacement relation of the classical von Mises truss (see \cite{Itskov-MoSM} for more details). In Sect. \ref{P-ratio} we apply
	the proposed constitutive equation to uniaxial tension and show that depending on the material constants it yields arbitrary values of Poisson's ratio except of $-1$. Other loading cases as for example simple shear and pure dilation are finally considered in Sect. \ref{other-lc}.

	\section{Model formulation}\label{Model}
	A general constitutive relation for a simple (Cauchy) elastic material can be given in terms of a generalized strain tensor $\mathbf E$ in the material description and the corresponding work-conjugate stress tensor $\mathbf T$ as 
	\begin{equation}\label{S(E)}
		\mathbf T = \hat{\mathbf T}\left( \mathbf E\right), 
	\end{equation}
	where $\hat{\mathbf T}$ denotes a response function.
	Recall that the tensors $\mathbf E$ and $\mathbf T$ do not change under a change of frame so that the relation of the form \eqref{S(E)} is a priori objective.
	
	Assuming the function $\hat{\mathbf T}\left( \mathbf E\right)$ to be analytic it can be developed into Taylor series around a stress and strain free reference configuration so that
	\begin{equation}\label{series-0f-eps}
		{\mathbf T}=\mathcal{E}_1 : \mathbf E +\mathcal{E}_2 : \mathbf E^2+\mathcal{E}_3 : \mathbf E^3 + \ldots, 
	\end{equation}
	where $\mathcal{E}_i, i=1,2\ldots$ represent fourth-order elasticity tensors and colon denotes a linear mapping of one second-order tensor into another one
	by a fourth-order tensor (see, e.g. \cite{Itskov}). Linearization of the series on the right hand side of \eqref{series-0f-eps} with respect to $\mathbf E$ would lead to ${\mathbf T}=\mathcal{E}_1 : \mathbf E$ which further reduces to the classical generalized Hooke law ${\tenf \sigma}=\mathcal{E}_1 : \tenf \epsilon$ at infinitesimal strains. For this reason, the linear term in $\mathbf E$ should be avoided. The next remaining term from the left is quadratic. It would
	lead to an unphysical response. Indeed, changing the sign of strain would not affect stresses while application of compressive (negative) stresses would not be possible
	at all. The lowest order reasonable nonlinear term is thus cubic. Since higher order
	terms in $\mathbf E$ would be negligible at infinitesimal strains we consider the following constitutive equation
	\begin{equation}\label{eps-power-3}
		\mathbf T=\mathcal{E}_3 : \mathbf E^3. 
	\end{equation}
	As a critical point we note that the above material model is characterized by zero stiffness including zero Young's and bulk modulus in the reference state. Thus, any 
	preload as for example by gravitation or air pressure would cause considerable deformation. The constitutive equation \eqref{eps-power-3} rewritten with respect to this new deformed reference configuration will again include a linear term (as in \eqref{series-0f-eps}) and consequently imply
	the superposition principle. For this reason, the strict form of the model \eqref{eps-power-3} would be possible for a solid in gravitation free space and vacuum. However, for light weight open-porous (meta-)materials (as for example aerogels) air pressure does not apply while the gravitation effect is insignificant so that the above mentioned linear term  can be neglected. As an example of a structure with zero stiffness
	in the reference state and a cubic force-displacement response similar to \eqref{eps-power-3} we can also mention the von Mises truss (see \cite{Itskov-MoSM} for more details) with zero original slope.

	Cauchy-elastic constitutive relations can lead to a path-dependent response and a release or consumption of energy in a closed
	loading loop which is critical from the thermodynamic point of view (see e.g. \cite{Carr}). To avoid this problem, we are going to formulate a hyperelastic (Green-elastic) material
	law leading to the cubic constitutive relation of the form \eqref{eps-power-3}. Indeed, a hyperelastic material law based on a positive-definite strain energy function
	would be path-independent and thermodynamically consistent. 
	
	A strain energy function leading to the constitutive relation of the form \eqref{eps-power-3} should be of order four in  $\mathbf E$. According to the classical invariant theory an isotropic scalar-valued tensor function of $\mathbf E$ can generally be expressed in terms of the principal
	traces $\text{tr}\mathbf E$, $\text{tr}\mathbf E^2$ and $\text{tr}\mathbf E^3$, where $\text{tr} \bullet$ denotes the trace of a second-order tensor. Thus, a fourth-order isotropic strain energy function of $\mathbf E$ 
	can generally be written by
	\begin{equation}\label{Psi-4}
		\Psi\left(\mathbf E  \right) = \alpha_1 \text{tr}\mathbf E^3 \text{tr}\mathbf E 
		+\alpha_2 \left( \text{tr}\mathbf E^2 \right)^2 
		+\alpha_3\text{tr}\mathbf E^2\left( \text{tr}\mathbf E\right)^2  
		+\alpha_4\left( \text{tr}\mathbf E\right)^4,
	\end{equation}
	where  $\alpha_i, i=1,2\ldots 4$ represent material constants. In the case of incompressible behavior or isochoric loading under infinitesimal strains $\text{tr}\mathbf E=\text{tr}\tenf{\epsilon} =0$, where $\tenf{\epsilon}$ denotes the Cauchy strain. Thus, $\Psi\left(\mathbf E\right) $ \eqref{Psi-4} reduces to the single second term, which requires that $\alpha_2>0$.
	In the following we consider a special form of \eqref{Psi-4} given by
	\begin{equation}\label{Psi-prop}
		\Psi\left(\mathbf E  \right) = \left[a \text{tr}\mathbf E^2 + b \left( \text{tr}\mathbf E\right)^2\right]^2, 
	\end{equation}
	where $a\ne0$ and $b$ represent material constants. Due to the quadratic form the strain energy \eqref{Psi-prop} is always positive-definite
	for all possible values of $a$ and $b$.  We also observe that \eqref{Psi-prop} represents a special case of \eqref{Psi-4} when $\alpha_1=0$, $\alpha_2=a^2$, $\alpha_3=2ab$ and $\alpha_4=b^2$.
	
	Alternatively, $\Psi\left(\mathbf E  \right)$ \eqref{Psi-prop}
	can be written in terms of the generalized principal strains $E_i, i=1,2,3$ which represent eigenvalues of 
	\begin{equation}\label{e-eigenvalues}
		\mathbf{E}=\left[
		\begin{array}{ccc}
			E_1 & 0 & 0 \\
			0 & E_2 & 0 \\
			0 & 0 & E_3
		\end{array}
		\right] \vect{e}_i\otimes \vect{e}_j,
	\end{equation}
	where $\vect{e}_i, i=1,2,3$ are orthonormal vectors and $\otimes$ denotes the tensor product of vectors.
	Accordingly,
	\begin{align}\label{Psi-princip}
		\Psi\left(E_1, E_2, E_3 \right)=
		&\left[a \left(E_1^2+ E_2^2+ E_3^2 \right)\right. \\ \nonumber
		&\left. + b \left(E_1+ E_2+ E_3 \right)^2\right]^2. 
	\end{align}
	
	A constitutive relation resulting from \eqref{Psi-prop} can be written as
	\begin{equation}\label{sigma-eps-gen}
		\mathbf{T}= \frac{\partial \Psi}{\partial {\mathbf E}} = 
		4\left(a {\mathbf E}+ b \text{tr}{\mathbf E} \,\mathbf{I} \right)
		\left[a \text{tr}{\mathbf E}^2 + b \left( \text{tr}{\mathbf E}\right)^2\right].
	\end{equation}
{\color{red}
In the following, we apply in \eqref{Psi-4} the logarithmic (Hencky) strain which can be obtained by
\begin{equation}\label{Hencky}
	\mathbf{E}=\frac{1}{2}\ln \mathbf{C}
\end{equation}
in terms of the right Cauchy-Green tensor $\mathbf{C}$. Eigenvalues of $\mathbf{E}$ \eqref{Hencky} can be expressed by $E_i=\ln \lambda_i$
in terms of the principal stretches $\lambda_i>0, i=1,2,3$. Accordingly, \eqref{Psi-prop} and \eqref{Psi-princip} result in
\begin{equation}\label{Psi-prop-log}
	\tilde{\Psi}\left(\mathbf C  \right) = \left[\frac{a}{4} \text{tr}\left( \ln\mathbf{C}\right)^2  + b  \ln^2 J\right]^2, 
\end{equation}
\begin{align}\label{Psi-princip-stretch}
	\hat{\Psi}\left(\lambda_1, \lambda_2, \lambda_3 \right)=
	&\left[a \left(\ln^2\lambda_1+ \ln^2\lambda_2+ \ln^2\lambda_3 \right)\right. \\ \nonumber
	&\left. + b  \ln^2 J\right]^2, 
\end{align}
where $J=\sqrt{\det \mathbf{C}}=\lambda_1\lambda_2\lambda_3$ denotes the relative volume change. Thus, due to the application of the logarithmic strain measure
\eqref{Hencky} the volumetric response of the model appears physically reasonable both at small and large strains. Indeed, according to \eqref{Psi-prop-log} and \eqref{Psi-princip-stretch} it tends to infinity when either $J\to 0$ or $J \to \infty$.
}

	\section{Poisson's ratio}\label{P-ratio}
	
	Poisson's ratio is defined in uniaxial tension as
	\begin{equation}\label{nu-def}
		\nu=\lim_{\epsilon_1\to 0}-\frac{\epsilon_2}{\epsilon_1},
	\end{equation}
	where $\epsilon_1$ denotes the Cauchy strain in the tension direction while for orthogonal directions $\epsilon_2=\epsilon_3$ in the case of isotropic material. At infinitesimal strains ${\mathbf E} \to \tenf{\epsilon}$ so that \eqref{nu-def} can alternatively be written in terms of the logarithmic strain as
	\begin{equation}\label{nu-def-E}
		\nu=\lim_{E_1\to 0}-\frac{E_2}{E_1}.
	\end{equation}
	Note that we deal here with the classical Poisson ratio defined at infinitesimal strains which is expressed in \eqref{nu-def} and \eqref{nu-def-E}
	by the limits at ${\epsilon_1\to 0}$ and ${E_1\to 0}$, respectively. In contrast, in nonlinear elasticity there are different measures of the lateral to longitudinal strain relation as for example apparent, secant, tangential Poisson's ratio or Poisson's function (see e.g. \cite{Miha} for more details). Moreover, contrary to $\nu$ \eqref{nu-def} there are no strict bounds for these measures at finite strains.
	
	By setting in \eqref{sigma-eps-gen} $E_2=E_3$ the normal stresses in the tension and lateral directions can be expressed, respectively, by
	\begin{align}\label{sigma-uni}
		T_1 =  4 &\left[ \left( a+b\right)E_1^2+4bE_1E_2+2\left( a+2b\right)E_2^2\right]\nonumber\\
		& \times\left[ \left( a+b\right)E_1+2bE_2\right],
	\end{align}
	\begin{align}
		T_2=T_3=  4 &\left[ \left( 2a+4b\right)E_2^2+4bE_1E_2+\left( a+b\right)E_1^2\right]\nonumber\\
		& \times\left[ \left( a+2b\right)E_2+bE_1\right].
	\end{align}
	{\color{red} Note that in the case of the logarithmic strain applied here $T_i, i=1,2,3$ represent the principal Kirchhoff stresses.}
	The condition of the stress-free lateral directions requires that $T_2=T_3=0$ and consequently
	\begin{align}\label{cubic-eq-nu}
		&\left[2\left(2\frac{b}{a}+1 \right) \left( \frac{E_2}{E_1}\right)^2 +4\frac{b}{a} \left( \frac{E_2}{E_1}\right) +\frac{b}{a}+1\right]\nonumber\\
		&\times\left[\left(2\frac{b}{a}+1 \right)\frac{E_2}{E_1}+\frac{b}{a}\right] =0.
	\end{align}
{\color{red}	Note that this equation defines the relation 
	\begin{equation}\label{r}
		r=-\frac{E_2}{E_1}=-\frac{\ln \lambda_2}{\ln\lambda_1}
	\end{equation}
 in terms of material constants independent of strains. Thus, the roots of \eqref{cubic-eq-nu} are valid both for small and large strains. At large strains $r$ is usually referred to as the biaxiality parameter since $\lambda_2=\lambda_1^{-r}$. At infinitesimal strains it collapses to Poisson's ratio $\nu=\lim_{E_1\to 0}r$.} Accordingly, with the definition \eqref{nu-def-E} and introducing a dimensionless material parameter $\beta=\frac{b}{a}$ we get 
	\begin{equation}\label{eqn-nu-uni}
		\left[\left(\nu -\frac{1}{2}\right)^{2} \beta +\frac{\nu^{2}}{2}+\frac{1}{4}\right]\left[\left(\nu -\frac{1}{2}\right) \beta +\frac{\nu}{2}
		\right] =0.
	\end{equation}
	The roots of this cubic equation can be given by
	\begin{equation}\label{roots-cubic}
		\nu_1=\frac{\beta}{1+2\beta}, \quad \nu_{2/3}=\frac{2 \beta \pm\sqrt{-6 \beta -2}}{2+4 \beta} .
	\end{equation}
	
	One can easily check that the expressions for $\nu_{2/3}$ \eqref{roots-cubic}$_2$ lead to the trivial result $S_1=0$ when inserted into
	\eqref{sigma-uni} as $E_2=-\nu_{2/3}E_1$. For this reason we focus in the following on the first root in \eqref{roots-cubic}.
	It results in the following expression of the stress in the tension direction:
	\begin{equation}\label{s1-uni}
		T_1=\frac{4 a^{2} \left(3 \beta +1\right)^{2}}{\left(1+2 \beta \right)^{2}}E_1^3=4a^2\eta E_1^3,
	\end{equation}
	where $\eta=4\frac{ \left(3 \beta +1\right)^{2}}{\left(1+2 \beta \right)^{2}}$ represents a stress factor. 
	In Fig. \ref{nu-beta} $\nu$ \eqref{roots-cubic}$_1$  and $\eta$ \eqref{s1-uni} are plotted versus $\beta$. Accordingly, the value $\beta=-\frac{1}{3}$ leads to $\eta=0$ and should thus be avoided. It corresponds to the value of Poisson's ratio $\nu=-1$ which is not possible in an isotropic material also for another reason.
	Indeed, let us consider an isotropic material with $\nu=-1$ under uniaxial tension. Due to the relations $E_2=E_3=-\nu E_1=E_1$ it causes pure dilatation which can appear in an isotropic material only under hydrostatic tension. However,
	in uniaxial tension one direction is subject to stress and two other lateral ones are stress-free, which violates isotropic symmetry of the material.
	Another special value of $\beta=-\frac{1}{2}$ should be avoided as well. In this case, the cubic equation \eqref{eqn-nu-uni} reduces to
	a quadratic one and the expression for the first root in \eqref{roots-cubic} is no more valid. Moreover, $\eta \to \infty$ when $\beta \to -\frac{1}{2} $ so that the material becomes infinitely stiff in uniaxial tension. 
	
	Apart from the exception of $-1$, any arbitrary values of $\nu$ can be reached according to \eqref{roots-cubic}$_1$ (see also Fig. \ref{nu-beta}). For example, within the interval $\beta\in[-\infty,-0.5[$ Poisson's ratio $\nu>0.5$. For example, $\nu=2$ for $\beta = -2/3$. This means that the volume decreases under uniaxial tension which is somewhat unusual. 
	The incompressibility limit with $\nu=0.5$ is reached when either $\beta\to\infty$ or $\beta\to-\infty$.  Indeed, $\lim_{\beta\to-\infty}\nu_1=0.5^+$ while $\lim_{\beta\to\infty}\nu_1=0.5^-$. Within the interval  
	$\beta\in]-\frac{1}{2},-\frac{1}{3}[$ Poisson's ratio $\nu<-1$. This is not possible for isotropic materials described by 
	a linearizable constitutive relation which reduces to 
	the classical generalized Hooke law at infinitesimal strains. For example, $\nu=-2$ for $\beta = -0.4$.
	\begin{figure}[!ht]
		\centering
		\includegraphics[width=.99\linewidth]{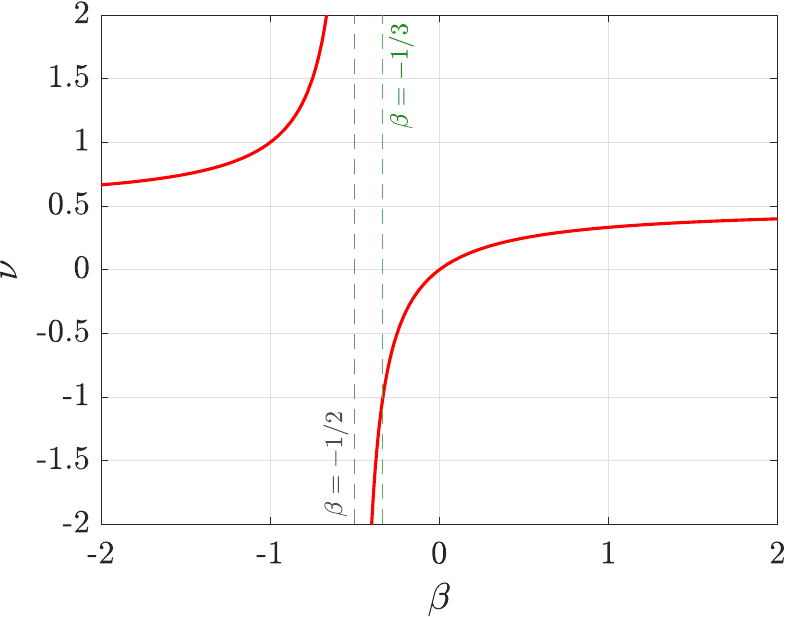} 
		\includegraphics[width=.99\linewidth]{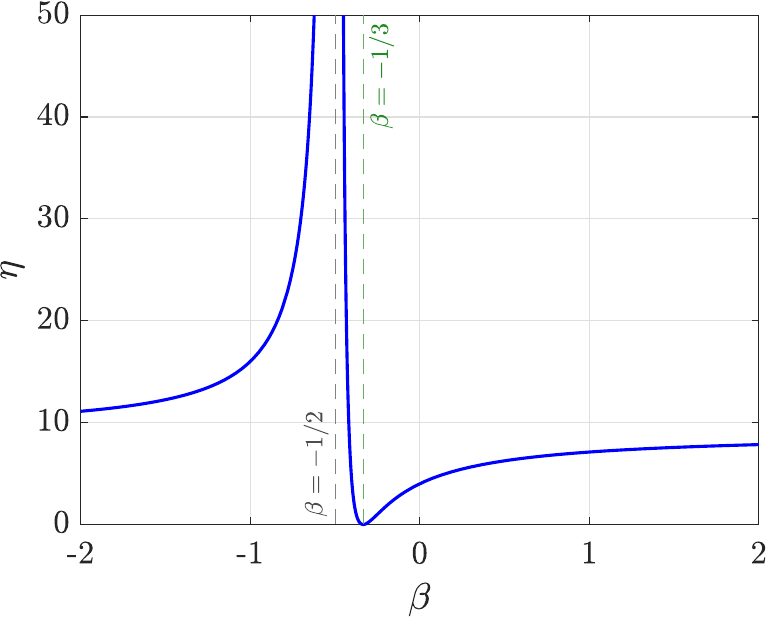} 
		\caption{Poisson's ratio (top) and the stress factor $\eta$ \eqref{s1-uni} (bottom) plotted versus the material parameter $\beta$ }
		\label{nu-beta}
	\end{figure}
	{\color{red}

The first Piola-Kirchhof stress resulting from \eqref{s1-uni} can be expressed by
\begin{equation}\label{P1-uni}
	P_1=4a^2\eta \frac{\ln^3\lambda_1}{\lambda_1}.
\end{equation}
One can easily see that the function $P_1\left( \lambda_1\right)$ is increasing within the interval $\lambda_1\in ]0,e^3]$ and begin to decrease only at very large stretches $\lambda_1>e^3\approx 20.1$. In Fig. \ref{P1-lambda} this function is plotted in normalized form versus $\lambda_1$.
	\begin{figure}[!ht]
	\centering
	\includegraphics[width=.99\linewidth]{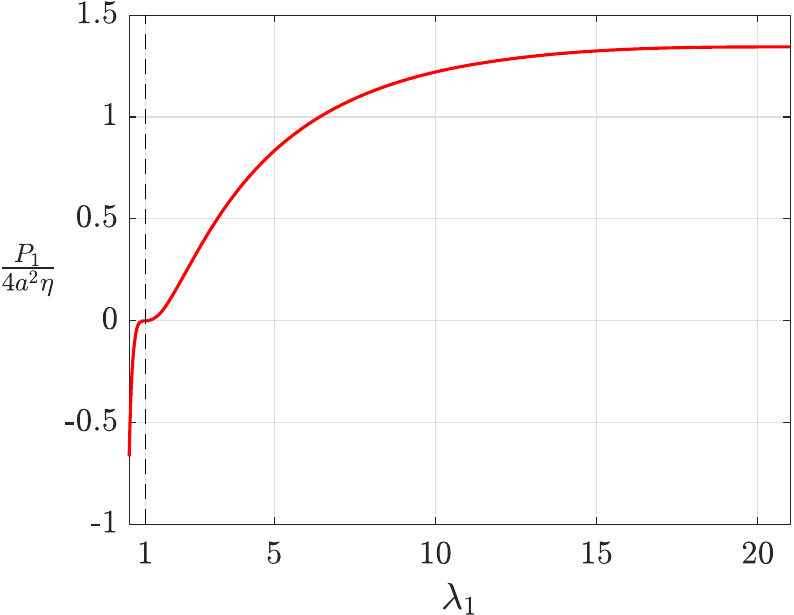} 
	\caption{First Piola-Kirchhoff stress versus stretch under uniaxial tension/compression}
	\label{P1-lambda}
\end{figure}

Finally, we study convexity properties of the strain energy function \eqref{Psi-princip-stretch}. The first necessary condition of rank-one convexity in Aubert \cite{Aubert} reads as
\begin{equation}
	\frac{\partial^2 \hat{\Psi}}{\partial \lambda_1^2}\ge 0\quad \forall \lambda_i>0.
\end{equation}
Evaluating the expression for $\lambda_2=\lambda_3=\lambda_1^{-r}$, where $r=\frac{\beta}{1+2\beta}$ we arrive at the condition
\begin{equation}
	-\left(3\beta + 1\right)^{2} \ln \lambda_1+6\beta^{3}+29 \beta^{2}+18 \beta +3 \ge 0,
\end{equation} 
which is not satisfied whenever $\ln\lambda_1>\frac{2 \beta^{2}+9 \beta +3}{3 \beta +1}$. Thus, the strain energy function \eqref{Psi-princip-stretch} is not rank-one convex and for this reason is not polyconvex.

}
	\section{Other deformation states}\label{other-lc}
	
	Further, we consider the response of the model in simple shear where {\color{red}the principal stretches are given by
	\begin{equation}\label{e-simple-shear}
	\lambda_{1/2}=\left( \frac{\sqrt{4+\gamma^2}\pm \gamma}{2}\right)^2, \quad \lambda_3=1 
	\end{equation}
	and $\gamma$ denotes the amount of shear. In this case, $J=1$ so that the strain energy function \eqref{Psi-princip-stretch} takes the form		
\begin{align}\label{Psi-shear}
	\hat{\Psi}\left(\lambda_1, \lambda_2, \lambda_3 \right)=
64a^2 {\ln^4\left(\frac{\sqrt{\gamma^{2}+4}}{2}+\frac{\gamma}{2}\right)}.
\end{align}
Thus, the Cauchy shear stress (work conjugate to $\gamma$) can be expressed by
	\begin{equation}\label{tau-shear}
		\tau=\frac{\partial \Psi}{\partial \gamma} = \frac{256a^2}{\sqrt{\gamma^{2}+4}} {\ln^3\left(\frac{\sqrt{\gamma^{2}+4}}{2}+\frac{\gamma}{2}\right)}.
	\end{equation}
This shear stress response is illustrated in Fig. \ref{tau-gamma}. One can easily see that $\tau$ is monotone increasing function of $\gamma$.
	\begin{figure}[!ht]
	\centering
    \includegraphics[width=.99\linewidth]{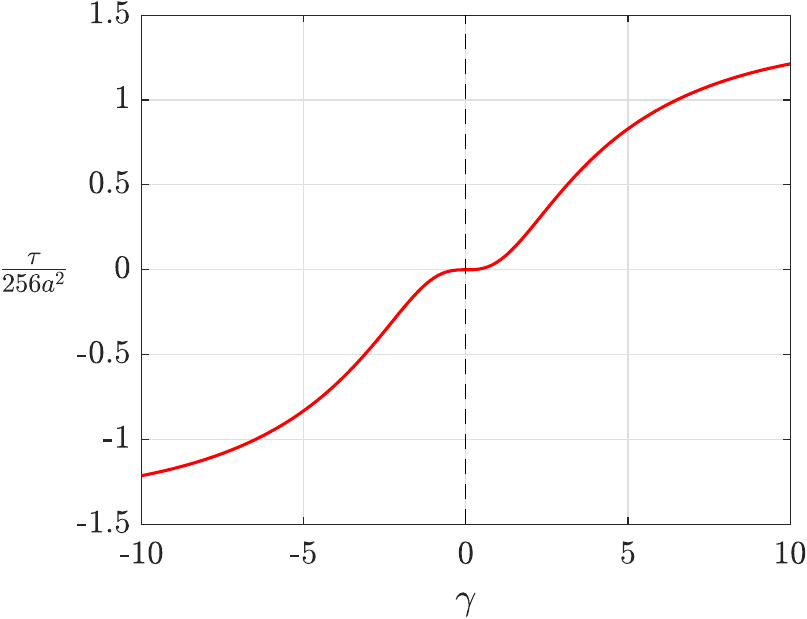} 
	\caption{Cauchy shear stress versus $\gamma$ under simple shear}
	\label{tau-gamma}
\end{figure}

Under equitriaxial tension/compression (pure dilatation) $E_1=E_2=E_3=\ln\lambda$ and consequently
	\begin{equation}
		\Psi= 9\left(a + 3b\right)^2 E^4,
	\end{equation}
	\begin{align}
		T_1=T_2=T_3= \frac{1}{3}\frac{\partial \Psi}{\partial E} &= 12\left(a + 3b\right)^2 E^3 \nonumber \\
		&=12\left(a + 3b\right)^2\ln^3\lambda.
	\end{align}
This stress response is of the same form as in uniaxial tension \eqref{s1-uni}. This concerns also the form \eqref{P1-uni} of the corresponding first Piola-Kirchhoff stress illustrated in Fig. \ref{P1-lambda} as well as of the Cauchy stress (hydrostatic tension/pressure under pure dilatation) as a function of $J$. Indeed, it can be expressed by
\begin{equation}
	-p=\frac{T_1}{J}=\frac{4}{9}\left(a + 3b\right)^2\frac{\ln^3J}{J},
\end{equation} 	
where $J=\lambda^3$.
}

	We can observe that the material response in all these loading cases is plausible and stable for all values of the material constants except of $\beta =-\frac{1}{3}$ and $\beta =-\frac{1}{2}$. This is in contrast to the classical linear (Lam\'{e}) elasticity where a Poisson's ratio above one half or less than minus one would lead to the negative bulk and shear modulus, respectively, and consequently to an unphysical response contradicting to the laws of thermodynamics. 
	
	\section{Conclusion}
	We proposed a hyperelastic isotropic material model described by a fourth-order strain energy function in terms of the {}\color{red}
	logarithmic (Hencky)} strain. This strain energy function represents a quadratic form and is thus a priori positive-definite. The model leads to a constitutive equation of third order in strain which cannot be linearized. For this reason, the superposition principle and as a result the classical Hooke law do not apply even at infinitesimal strains. 
	Another unusual consequence of the proposed model is that Poisson's ratio can take arbitrary values except of $-1$. Thus, it can be smaller than
	$-1$ or greater than $0.5$ in full agreement with the laws of thermodynamics. For all these values responses of the material model in uniaxial tension, simple shear and pure dilation appear reasonable and stable. {\color{red}Nevertheless, under uniaxial tension the strain energy function is shown not to be rank-one convex and for this reason it is not polyconvex. }
	As another disadvantage of the model we note zero stiffness including zero Young's and bulk modulus in the reference state. Thus, any preload as for example by gravitation or air pressure would cause considerable deformation and induce again a linear term into the constitutive equation. Thus, the proposed material model could be useful for light weight open-porous materials or metamaterials.
	
	\section*{Acknowledgment}
	
	The author would like to thank Rajesh Chandrasekaran for his assistance with the preparation of figures for this paper. {\color{red} The author is also very grateful to anonymous reviewers for valuable comments to the previous version of the paper.} 
	
\section*{Declarations}

The author declares that he has no known competing financial interests or personal relationships that could have appeared to influence the work reported in this paper.
\\ \\
Data Availability declaration: there is no data supporting this study.\\

\bibliographystyle{cas-model2-names}



\end{document}